\documentclass[smallextended,natbib,runningheads]{svjour3}
\smartqed  
\usepackage{graphicx}
\usepackage{amssymb}
\usepackage{amsmath}
\usepackage{amsfonts}

\begin{document}

\title{Asymptotic theory   of multiple-set linear  canonical analysis}


\author{ Guy Martial Nkiet}


\institute{
            Guy Martial Nkiet \at
 Universit\'e des Sciences et Techniques de Masuku, D\'epartement de Math\'ematiques et Informatique, BP 943 Franceville, Gabon.
            \email{gnkiet@hotmail.com}
}

\date{Received: date / Revised: date}

\maketitle

\begin{abstract}
 This paper deals  with asymptotics for multiple-set linear canonical  analysis (MSLCA). A definition of this analysis, that adapts the classical one to the context of Euclidean random variables, is given and properties of the related canonical coefficients are derived. Then, estimators of the MSLCA's elements, based on empirical covariance operators,  are proposed and asymptotics for these estimators are obtained. More precisely, we prove their consistency and we obtain asymptotic normality for the estimator of the operator that gives MSLCA, and also for the estimator of the vector of canonical coefficients. These results are then used to obtain a test for mutual non-correlation between the involved Euclidean random variables. 
\end{abstract}

\section{Introduction}

Multiple-set linear canonical analysis (MSLCA), also known as generalized canonical correlation  analysis, has been extensively discussed in the literature, see Kettenring (1971),  Gifi (1991), Gardner et al. (2006), Takane et al. (2008), Tenenhaus and Tenenhaus (2011), as well as the further references contained therein. It is a statistical method that generalizes linear canonical analysis (LCA) to the case where more than two sets of variables are considered, which is of  a real interest since in applied statistical studies it is common to collect data from the observation of several sets of variables on a given population. However, although this interest, several aspects under which LCA  has been studied have never  been addressed to MSLCA. For example, asymptotic theory for LCA and related applications  have been tackled by several authors (e.g., Muirhead and Waternaux (1980), Anderson (1999), Pousse (1992), Fine (2000), Dauxois et al. (2004)). It would be natural to wonder how the obtained results extend to the case of MSLCA but, to the best of our knowledge,  such an approach has  never been tackled.

In this paper, we introduce an asymptotic theory for MSLCA. For doing that, we first define  in Section 2 the notion of MSLCA for Euclidean random variables, that is random variables valued into  Euclidean vector spaces. This analysis is defined from  a maximization problem under specified constraints, and  shown to be  obtained from spectral analysis of a suitable operator. Properties of the related eigenvalues, called canonical coefficients, are then given. In Section 3, we tackle the problem of estimating MSLCA. More precisely, estimators based on empirical covariance operators are introduced. Then, consistency of the obtained estimators is proved. Further, we derive the asymptotic distribution of the used estimator of the aforementioned operator, and also that of  the estimator of the vector of canonical coefficients in the general case as well as in the case of elliptical distribution. Section 4 is devoted to the introduction of a test for mutual non-correlation between the random variables involved in MSLCA. The results obtained for asymptotic theory of MSLCA are then used in order to derive the asymptotic distribution of the used test statistic under null hypothesis.

\section{Multiple-set canonical linear analysis of Euclidean random variables}
\noindent For an integer $K\geq 2$, let us  consider random variables $X_1,\cdots,X_K$ defined on a probability space $(\Omega,\mathcal{A},P) $ and valued into  Euclidean vector spaces $\mathcal{X}_1,\cdots,\mathcal{X}_K$ respectively. Denoting by $\mathbb{E}$ the mathematical expectation related to $P$, we assume that, for any $k\in\{1,\cdots, K\}$, we have $\mathbb{E}(\Vert X_k\Vert _k^2)<+\infty $ where $\Vert\cdot\Vert _k$ denotes the norm induced by the inner product  $\langle\cdot ,\cdot\rangle _ k$ of $\mathcal{X}_k$, and, without loss of generality, that   $\mathbb{E}(X_k)=0$.  Each vector $\alpha$ in the vector space $\mathcal{X}:=\mathcal{X}_1\times\cdots\mathcal{X}_K$ will be  writen as
\[
\alpha=\left(
\begin{array}[c]{c}
\alpha_1\\
\vdots\\
\alpha_K
\end{array}
\right),
\]
 and we recall that $\mathcal{X}$ is an Euclidean  vector space equipped with the inner product $\langle\cdot,\cdot\rangle_\mathcal{X}$ defined by:
\[
\forall\alpha\in \mathcal{X},\,\forall\beta\in\mathcal{X},\,\,\langle\alpha,\beta\rangle_\mathcal{X}=\sum_{k=1}^K\langle\alpha_k,\beta_k\rangle_k.
\]
We denote by  $\Vert\cdot\Vert _\mathcal{X}$  the norm induced by this inner product. Considering the $\mathcal{X}$-valued random variable 
\[
X=\left(
\begin{array}[c]{c}
X_1\\
\vdots\\
X_K
\end{array}
\right),
\]
we can give the following definition which adapts the classical definition of multiple-set canonical analysis (e.g., Gifi (1991), Gardner et al. (2006), Takane et al. (2008)) to the context of Euclidean random variables.

\bigskip

\noindent\textbf{Definition 2.1.}\textsl{ The multiple-set linear canonical analysis
(MSLCA) of  $X$ is the search of a sequence  $\left(\alpha^{(j)}\right)_{1\leq j\leq q}$ of vectors of  $E$, where $q=\dim(\mathcal{X})$,  satisfying:
\begin{equation}\label{sol1}
\alpha^{(j)}=\arg\max_{\alpha\in C_j}\mathbb{E}\left(<\alpha , X>_ \mathcal{X}^2\right),
\end{equation}
where
\begin{equation}\label{contr1}
C_1=\left\{\alpha\in \mathcal{X}\,/\,\sum_{k=1}^Kvar(<\alpha_k,X_k>_k)=1\right\},
\end{equation}
and, for  $j\geq 2$ :
\begin{equation}\label{contr2}
C_j=\left\{\alpha\in  C_1\,/\,\sum_{k=1}^Kcov\left(<\alpha^{(r)}_k,X_k>_k,<\alpha_k,X_k>_k\right)=0,\,\,\forall r\in\{1,\cdots,j-1\}\right\}.
\end{equation}}

\bigskip

\noindent\textbf{Remark 2.1}

\noindent 1) The constraints sets given in (\ref{contr1}) and (\ref{contr2}) can be expressed by using covariance operators defined for $(k,\ell)\in\{1,\cdots,K\}^2$  by:
\[
V_{k\ell}=\mathbb{E}\left(X_\ell\otimes X_k\right)=V_{\ell k}^\ast \textrm{ and } V_k:=V_{kk},
\]
where   $\otimes $ denotes the tensor product such that, for any $(x,y)$, $x\otimes y $ is the linear map : $h\mapsto <x,h>y$, and $T^\ast $ denotes the adjoint of $T$. Indeed, it is easily seen that, for $(\alpha,\beta)\in\mathcal{X}^2$,
\[
var\left(\langle\alpha_k,X_k\rangle_k\right)=
\mathbb{E}\left(\langle\alpha_k,X_k\rangle_k^2\right)
=\mathbb{E}\left(\langle\alpha_k,(X_k\otimes X_k)(\alpha_k)\rangle_k\right)
=\langle\alpha_k,V_k\alpha_k\rangle_k,
\]
and
\begin{eqnarray*}
cov\left(\langle\alpha_k,X_k\rangle_k,\langle\beta_\ell,X_\ell\rangle_\ell\right)&=&
\mathbb{E}\left(\langle\alpha_k,X_k\rangle_k\langle\beta_\ell,X_\ell\rangle_\ell\right)
=\mathbb{E}\left(\langle\alpha_k,(X_\ell\otimes X_k)(\beta_\ell)\rangle_k\right)\\
&=&\langle\alpha_k,V_{k\ell}\beta_\ell\rangle_k.
\end{eqnarray*}
Therefore,
\begin{equation}\label{contr12}
C_1=\left\{\alpha\in \mathcal{X}\,/\,\sum_{k=1}^K<\alpha_k,V_k\alpha_k>_k=1\right\},
\end{equation}
and
\begin{equation}\label{contr22}
C_j=\left\{\alpha\in  C_1\,/\,\sum_{k=1}^K<\alpha^{(r)}_k,V_k\alpha_k>_k=0,\,\,\forall r\in\{1,\cdots,j-1\}\right\}.
\end{equation}

\noindent 2) For any $\alpha\in C_1$, one has:
\begin{eqnarray*}
\mathbb{E}\left(<\alpha , X>_\mathcal{X}^2\right) &=&\mathbb{E}\left(\left(\sum_{k=1}^K<\alpha_k , X_k>_k\right)^2\right)
=\sum_{k=1}^K\sum_{\ell=1}^K\mathbb{E}\left(<\alpha_k , X_k>_k<\alpha_\ell , X_\ell>_\ell\right)\\
&=&\sum_{k=1}^K\mathbb{E}\left(<\alpha_k , X_k>_k^2\right)+\sum_{k=1}^K\sum_{\stackrel{\ell=1}{\ell\neq k}}^K\mathbb{E}\left(<\alpha_k , X_k>_k<\alpha_\ell , X_\ell>_\ell\right)\\
&=&\sum_{k=1}^Kvar(<\alpha_k,X_k>_k)+\sum_{k=1}^K\sum_{\stackrel{\ell=1}{\ell\neq k}}^K<\alpha_k,V_{k\ell}\alpha_\ell>_k\\
&=&1+\sum_{k=1}^K\sum_{\stackrel{\ell=1}{\ell\neq k}}^K<\alpha_k,V_{k\ell}\alpha_\ell>_k=1+\varphi(\alpha),
\end{eqnarray*}
where 
\begin{equation}\label{phi}
\varphi(\alpha)=\sum_{k=1}^K\sum_{\stackrel{\ell=1}{\ell\neq k}}^K<\alpha_k,V_{k\ell}\alpha_\ell>_k.
\end{equation}
Then, the MSLCA of $X$ is obtained by minimizing  $\varphi(\alpha)$ under the constraints expressed in  (\ref{contr12}) and (\ref{contr22}). 

\bigskip

\noindent For $k\in\{1,\cdots,K\}$, the covariance operator $V_k$ is a self-adjoint non-negative operator. From now on, we assume that it is invertible. Let $\tau_k$ be the canonical projection defined as 
\[
\tau_k\,:\,\alpha\in \mathcal{X}\mapsto \alpha_k\in \mathcal{X}_k;
\]
its  adjoint $\tau_k^\ast$ of $\tau_k$ is the map  given by:
\[
\tau_k^\ast\,:\,t\in \mathcal{X}_k\mapsto (\underbrace{0,\cdots,0}_{k-1\textrm{ times}},t,0,\cdots , 0)^T\in \mathcal{X},
\]
where we denote by $a^T$ the transposed of  $a$. Now,
 let us consider the operators of $\mathcal{L}(\mathcal{X})$  given by:
\[
\Phi=\sum_{k=1}^K\tau_k^\ast V_k\tau_k\,\,\,\,
\textrm{ and }\,\,\,\,
\Psi=\sum_{k=1}^K\sum_{\stackrel{\ell=1}{\ell\neq k}}^K\tau_k^\ast V_{k\ell}\tau_\ell.
\]
From the fact that $\tau_k\tau^\ast_ \ell=\delta_{k\ell}I_k$, where $\delta_{k\ell}$ is the usual Kronecker symbol and $I_k$ is the identity operator of $\mathcal{X}_k$, it  is easily seen that  $\Phi$ is also an invertible self-adjoint and non-negative operator, with $\Phi^{-1}=\sum_{k=1}^K\tau_k^\ast V_k^{-1}\tau_k$ and $\Phi^{-1/2}=\sum_{k=1}^K\tau_k^\ast V_k^{-1/2}\tau_k$. The following theorem   shows how to obtain a MSLCA of $X$.  It just repeats a known result (e.g., Gifi (1991), Takane et al. (2008)) within the framework  used for  this paper.
\bigskip

\noindent\textbf{Theorem 2.1. } \textsl{Letting $\left\{
\beta^{(1)},\cdots,\beta^{(q)}\right\} $  be an orthonormal basis of   $\mathcal{X}$  such that  $\beta^{(j)}$ is an eigenvector of the operator $T=\Phi^{-1/2}\Psi\Phi^{-1/2}$ associated with the  $j$-th largest eigenvalue  $\rho_j$ of $T$. Then, the sequence  $\left(\alpha^{(j)}\right)_{1\leq j\leq q}$ given by: 
\[
\alpha^{(j)}=\Phi^{-1/2}\beta^{(j)}=\left(V_{1}^{-1/2}\beta^{(j)}_1,\cdots,V_{K}^{-1/2}\beta^{(j)}_K\right),
\]
consists of  solutions of (\ref{sol1}) under the constraints (\ref{contr1}) and (\ref{contr2}), and we have: $\rho_j=<\beta^{(j)},T\beta^{(j)}>_E=\varphi(\alpha^{(j)})$.}

\noindent\textit{Proof.}  Putting  $\beta_k=V_k^{1/2}\alpha_k$ and $\beta^{(r)}=V_k^{1/2}\alpha^{(r)}_k$, we have:
\begin{eqnarray}\label{fonct}
\varphi(\alpha)&=&\sum_{k=1}^K\sum_{\stackrel{\ell=1}{\ell\neq k}}^K<V_k^{-1/2}\beta_k,V_{k\ell}V_\ell^{-1/2}\beta_\ell>_k\nonumber\\
&=&\sum_{k=1}^K\sum_{\stackrel{\ell=1}{\ell\neq k}}^K<\beta_k,V_k^{-1/2}V_{k\ell}V_\ell^{-1/2}\beta_\ell>_k=:\psi(\beta),
\end{eqnarray}
where 
\[
\beta=\left(
\begin{array}{c}
\beta_1\\
\vdots\\
\beta_K
\end{array}
\right)\in \mathcal{X}. 
\]
Since  $V_k=V_k^{1/2}V_k^{1/2}$, having $\alpha\in C_j$ is equivalent to having  $\beta\in C_j^\prime$, where:
\begin{eqnarray}\label{contr13}
C_1^\prime
&=&\left\{\beta\in \mathcal{X}\,/\,\sum_{k=1}^K\Vert \beta_k\Vert_k^2=1\right\}=\left\{\beta\in \mathcal{X}\,/\,\Vert \beta\Vert_\mathcal{X}^2=1\right\},
\end{eqnarray}
and for $j\geq 2$:
\begin{eqnarray}\label{contr23}
C_j^\prime&=&\left\{\beta\in  C_1\,/\,\sum_{k=1}^K<\beta^{(r)}_k,\beta_k>_k=0,\,\,\forall r\in\{1,\cdots,j-1\}\right\}\nonumber\\
&=&\left\{\beta\in  C_1\,/\,<\beta^{(r)},\beta>_\mathcal{X}=0,\,\,\forall r\in\{1,\cdots,j-1\}\right\}.
\end{eqnarray}
Further, for any   $\beta\in\mathcal{X}$:
\begin{eqnarray*}
\Psi\Phi^{-1/2}\beta&=&\sum_{k=1}^K\sum_{\stackrel{\ell=1}{\ell\neq k}}^K\sum_{j=1}^K\tau_k^\ast V_{k\ell}\tau_\ell\tau_j^\ast V_j^{-1/2}\tau_j\beta
=\sum_{k=1}^K\sum_{\stackrel{\ell=1}{\ell\neq k}}^K\sum_{j=1}^K\delta_{\ell j}\tau_k^\ast V_{k\ell} V_j^{-1/2}\tau_j\beta\\
&=&\sum_{k=1}^K\sum_{\stackrel{\ell=1}{\ell\neq k}}^K\tau_k^\ast V_{k\ell} V_\ell^{-1/2}\tau_\ell\beta,
\end{eqnarray*}
and
\begin{eqnarray*}
\Phi^{-1/2}\Psi\Phi^{-1/2}\beta
&=&\sum_{k=1}^K\sum_{\stackrel{\ell=1}{\ell\neq k}}^K\sum_{j=1}^K\tau_j^\ast V_j^{-1/2}\tau_j\tau_k^\ast V_{k\ell} V_\ell^{-1/2}\tau_\ell\beta\\
&=&\sum_{k=1}^K\sum_{\stackrel{\ell=1}{\ell\neq k}}^K\sum_{j=1}^K\delta_{jk}\tau_j^\ast V_j^{-1/2} V_{k\ell} V_\ell^{-1/2}\tau_\ell\beta\\
&=&\sum_{k=1}^K\sum_{\stackrel{\ell=1}{\ell\neq k}}^K\tau_k^\ast V_k^{-1/2}V_{k\ell} V_\ell^{-1/2}\tau_\ell\beta\\
&=&\sum_{k=1}^K\sum_{\stackrel{\ell=1}{\ell\neq k}}^K\tau_k^\ast V_k^{-1/2}V_{k\ell} V_\ell^{-1/2}\beta_\ell .
\end{eqnarray*}

Thus,
\begin{eqnarray*}
<\beta,\Phi^{-1/2}\Psi\Phi^{-1/2}\beta >_\mathcal{X}&=&\sum_{k=1}^K\sum_{\stackrel{\ell=1}{\ell\neq k}}^K<\beta,\tau_k^\ast V^{-1/2}_kV_{k\ell}V_\ell^{-1/2}\beta_\ell>_\mathcal{X}\\
&=&\sum_{k=1}^K\sum_{\stackrel{\ell=1}{\ell\neq k}}^K<\tau_k\beta, V^{-1/2}_kV_{k\ell}V_\ell^{-1/2}\beta_\ell>_k\\
&=&\sum_{k=1}^K\sum_{\stackrel{\ell=1}{\ell\neq k}}^K<\beta_k, V^{-1/2}_kV_{k\ell}V_\ell^{-1/2}\beta_\ell>_k=\psi(\beta),
\end{eqnarray*}
where $\psi$ is defined in (\ref{fonct}). Then, the MSLCA optimization problem reduces to the maximization of   $<\beta,\Phi^{-1/2}\Psi\Phi^{-1/2}\beta >_\mathcal{X}$ under the constraints (\ref{contr13})  and  (\ref{contr23}). Since $T=\Phi^{-1/2}\Psi\Phi^{-1/2}$ is a self-adjoint operator, this is a well known maximization problem for which  a solution is obtained from the spectral analysis of $T$ as stated in the theorem. \hfill$\Box$

\bigskip

\bigskip

\noindent\textbf{Definition 2.2. }\textsl{The  $\rho_j$'s are termed the canonical coefficients. The $\alpha^{(j)}$'s are termed  vectors  of canonical directions}.

\bigskip

\noindent The following theorem gives some properties of the canonical coefficients.

\bigskip

\noindent\textbf{Theorem 2.2. }

\noindent\textsl{(i)} $\forall j\in\{1,\cdots,q\}$, $-1\leq\rho_j\leq K(K-1).$

\noindent\textsl{(ii)} $\forall j\in\{1,\cdots,q\}$, $\rho_j=0\Leftrightarrow\forall (k,\ell)\in\{1,\cdots,K\}^2,\,k\neq\ell,\,V_{k\ell}=0.$

\noindent\textit{Proof. }

\noindent (i) First,   using  (\ref{phi}), we have for any $ j\in\{1,\cdots,q\}$, 
\[
\rho_j=\varphi(\alpha^{(j)})=\mathbb{E}\left(<\alpha^{(j)},X>_\mathcal{X}^2\right)-1\geq-1.
\]
On the other hand, we have:
\begin{eqnarray*}
\rho_j=\varphi(\alpha^{(j)})&=&\sum_{k=1}^K\sum_{\stackrel{\ell=1}{\ell\neq k}}^K\mathbb{E}\left(<\alpha^{(j)}_k , X_k>_k<\alpha^{(j)}_\ell , X_\ell>_\ell\right)\\
&\leq &\sum_{k=1}^K\sum_{\stackrel{\ell=1}{\ell\neq k}}^K\sqrt{\mathbb{E}\left(<\alpha^{(j)}_k , X_k>^2_k\right)}\sqrt{\mathbb{E}\left(<\alpha^{(j)}_\ell , X_\ell>^2_\ell\right)}.
\end{eqnarray*}
Since, for any $k\in\{1,\cdots,K\}$, one has:
\[
\mathbb{E}\left(<\alpha^{(j)}_k , X_k>^2_k\right)=var\left(<\alpha^{(j)}_k , X_k>_k\right)\leq\sum_{\ell=1}^Kvar\left(<\alpha^{(j)}_\ell , X_\ell>_\ell\right)=1,
\]
it follows that:
\[
\rho_j\leq\sum_{k=1}^K\sum_{\stackrel{\ell=1}{\ell\neq k}}^K1=K(K-1).
\]
(ii) Since the  $\rho_j$'s are the eigenvalues of   $T$, we have:
\[
\forall j\in\{1,\cdots,q\},\, \rho_j=0\Leftrightarrow T=0\Leftrightarrow  \Psi=0\Leftrightarrow \forall (k,\ell)\in\{1,\cdots,K\}^2,\,k\neq\ell,\,V_{k\ell}=0.
\]
\hfill $\Box$

\bigskip

\noindent\textbf{Remark 2.2. } 

\noindent 1) When $K=2$, one has  $\Phi=\tau_1^\ast V_1\tau_1+\tau_2^\ast V_2\tau_2$ and  $\Psi=\tau_1^\ast V_{12}\tau_2+\tau_2^\ast V_{21}\tau_1$. Then it is  easy to check that $T=\tau_1^\ast S\tau_2+ \tau_2^\ast S^\ast \tau_1$, where $S=V_1^{-1/2}V_{12}V_2^{-1/2}$. Let $x$ be an eigenvector of $T$ associated with an eigenvalue $\rho\neq 0$. We have $Tx=\rho x$, that is equivalent to having:
\[
\tau_1^\ast\left(S\tau_2x-\rho \tau_1x\right)=-\tau_2^\ast\left(S^\ast\tau_1x-\rho \tau_2x\right).
\]
This implies:
\[
\left\{
\begin{array}{c}
S\tau_2x=\rho\, \tau_1x\\
S^\ast\tau_1x=\rho\, \tau_2x
\end{array}
\right.
\]
and, putting $x_1=\tau_1x$ and $x_2=\tau_2x$, we obtain 
\begin{equation}\label{lca}
x_2=\rho^{-1}S^\ast x_1\,\,\,\textrm{ and }\,\,\,Rx_1=\rho^2x_1,
\end{equation}
 where 
\[
R=SS^\ast=V_1^{-1/2}V_{12}V_2^{-1}V_{21}V_2^{-1/2}.
\]
Conversely, if  (\ref{lca}) holds  then,  putting $x=\tau_1^\ast x_1+\tau_2^\ast x_2$, we have:
\begin{eqnarray*}
Tx&=&\tau_1^\ast S\tau_2x+\tau_2^\ast S^\ast\tau_1x=\tau_1^\ast S x_2+\tau_2^\ast S^\ast x_1
=
\rho^{-1}\tau_1^\ast SS^\ast x_1+\rho\,\tau_2^\ast x_2\\
&=&\rho^{-1}\tau_1^\ast R x_1+\rho\,\tau_2^\ast x_2
=\rho\left(\tau_1^\ast x_1+\tau_2^\ast x_2\right)=\rho\,x.
\end{eqnarray*}
Moreover, since
\[
\Vert x_2\Vert_2=\rho^{-1}\Vert S^\ast x_1\Vert_1=\rho^{-1}\sqrt{< S^\ast x_1, S^\ast x_1>_2}=\rho^{-1}\sqrt{<S S^\ast x_1,  x_1>_1}=\Vert x_1\Vert_1
\]
and
\[
\Vert x\Vert_\mathcal{X}^2=\Vert x_1\Vert_1^2+\Vert x_2\Vert_2^2
\]
it follows that
\[
\Vert x_1\Vert_1=\Vert x_2\Vert_2=\frac{1}{\sqrt{2}}\Vert x\Vert_\mathcal{X}.
\]

\noindent 2) The preceding remark shows the equivalence between MSLCA and linear canonical analysis (LCA) when $K=2$. Recall that LCA of $X_1$ and $X_2$ is obtained from the spectral analysis of $R$ (see, e.g., Dauxois and Pouse (1975), Pousse (1992), Fine (2000)). More precisely,   $\left\{
\beta^{(j)},\rho_j\right\}_{1\leq j \leq q} $ is defined as in Theorem 2.1 if, and only if,  $\left\{
u_1^{(j)}, u_2^{(j)}\rho_j^2\right\}_{1\leq j \leq q}$, where $u_\ell^{(j)}=\frac{1}{\sqrt{2}}\tau_\ell \beta^{(j)}$ ($\ell\in\{1,2\}$), is a LCA of $X_1$ and $X_2$.

\section{Estimation and asymptotic theory}
In this section, we deal with estimation of MSLCA. For $k=1,\cdots, K$, let $\{X_k^{(i)}\}_{1\leq i\leq n}$  be  an i.i.d. sample  of $X_k$. We use empirical covariance operators  for defining estimators of MSLCA elements. Then, consistency and asymptotic normality are obtained for the resulting estimators of the vectors of canonical directions and the canonical coefficients.
\subsection{Estimation and almost sure convergence}
\noindent For $(k,\ell)\in\{1,\cdots,K\}^2$, let us consider the sample means and covariance operators:
\[
\overline{X}_{k\cdot n}=\frac{1}{n}\sum_{i=1}^nX^{(i)}_k,\,\,\,
\widehat{V}_{k\ell\cdot n}=\frac{1}{n}\sum_{i=1}^n\left(X^{(i)}_\ell-\overline{X}_{\ell\cdot n}\right)\otimes \left(X^{(i)}_k-\overline{X}_{k\cdot n}\right),\,\,\,
\widehat{V}_{k\cdot n}:=\widehat{V}_{kk\cdot n},
\]
and the random operators valued into $\mathcal{L}(\mathcal{X})$ defined as
\[
\widehat{\Phi}_n=\sum_{k=1}^K\tau_k^\ast \widehat{V}_{k\cdot n}\tau_k\,\,\,\,
\textrm{ and }\,\,\,\,
\widehat{\Psi}_n=\sum_{k=1}^K\sum_{\stackrel{\ell=1}{\ell\neq k}}^K\tau_k^\ast \widehat{V}_{k\ell\cdot n}\tau_\ell.
\]
Then, we estimate $T$ by  
\[
\widehat{T}_n=\widehat{\Phi}_n^{-1/2}\widehat{\Psi}_n\widehat{\Phi}_n^{-1/2}.
\] 
Considering the eigenvalues $\widehat{\rho}_{1\cdot n}\geq\widehat{\rho}_{2\cdot n}\cdots\geq\widehat{\rho}_{q\cdot n}$ of  $\widehat{T}_n$, 
and  $\left\{
\widehat{\beta}^{(1)}_n,\cdots,\widehat{\beta}^{(q)}_n\right\} $   
an orthonormal basis of   $\mathcal{X}$  such that  $\widehat{\beta}^{(j)}_n$ is an eigenvector of $\widehat{T}_n $ associated with $\widehat{\rho}_{j\cdot n}$. Then, we estimate $\rho_j$ by $\widehat{\rho}_{j\cdot n}$, and $\beta^{(j)}$ by $ \widehat{\beta}^{(j)}_n$. The following theorem establishes strong consistency for these estimators.

\bigskip

\noindent\textbf{Theorem 3.1. } \textsl{For any integer $j\in\{1,\cdots,q\}$ :}

\noindent\textsl{ (i)  $\widehat{\rho}_{j\cdot n}$ converge almost surely, as $n\rightarrow +\infty $, to
$\rho_j$.}

\noindent\textsl{ (ii)  $ sign(\langle\widehat{\beta}^{(j)}_n,\beta^{(j)}\rangle _\mathcal{X})\,\widehat{\beta}^{(j)}_n$  converges almost surely, as $n\rightarrow +\infty $, to $\beta^{(j)}$  in $\mathcal{X}$.}

\noindent\textit{Proof}. From obvious applications of the strong law of large numbers, it is easily seen that $\widehat{T}_n$  converges almost surely in $\mathcal{L}(\mathcal{X})$, as $n\rightarrow +\infty$ to $T$. Then using Lemma 1 in  Ferr\'{e} and Yao (2003), we obtain the inequality
$\vert  \widehat{\rho}_{j\cdot n}-\rho_j\vert  \leq\Vert  \widehat{T}_n-T\Vert$ 
from what $(i)$ is deduced. Clearly, each $\beta^{(j)}\otimes \beta^{(j)}$ and  is a projector onto an eigenspace. Thefore, using Proposition 3 in Dossou-Gbete and Pousse (1991), we deduce  that  $\widehat{\beta}^{(j)}_n
\otimes\widehat{\beta}^{(j)}_n$  converges almost surely in ${\mathcal{L}(\mathcal{X})}$  to 
$\beta^{(j)}\otimes\beta^{(j)}$, as $n\rightarrow +\infty$.  Using again Lemma 1 in  Ferr\'{e} and Yao (2003), we obtain the inequality
\[
\left\|  sign(\langle\widehat{\beta}^{(j)}_n,\beta^{(j)}\rangle _\mathcal{X})\,\widehat{\beta}^{(j)}_n-\beta^{(j)}\right\|_\mathcal{X} \leq2\sqrt{2}\left\| \widehat{\beta}^{(j)}_n
\otimes\widehat{\beta}^{(j)}_n-\beta^{(j)}\otimes \beta^{(j)}\right\| 
\]
from what we deduce  \textit{(ii)}.
\hfill $\Box$
\subsection{Asymptotic distribution}
\noindent  In this section, we assume that, for $k\in\{1,\cdots,K\}$, we have $\mathbb{E}\left(\Vert X_k\Vert_k^4\right)<+\infty$ and $V_k=I_k$, where $I_k$ denotes the identity operator of $\mathcal{X}_k$. We first derive an asymptotic distribution for   $\widehat{T}_n $,  then we  obtain these of the canonical coefficients.

\bigskip

\noindent\textbf{Theorem 3.2. }\textsl{  $\sqrt{n}\left(\widehat{T}_n -T\right)$ converges in distribution, as $n\rightarrow +\infty$, to a random variable $U$ having a normal distribution in $\mathcal{L}(\mathcal{X})$, with mean $0$ and covariance operator $\Gamma$ equal to that of the random operator:}
\[
Z=\sum_{k=1}^K\sum_{\stackrel{\ell=1}{\ell\neq k}}^K-\frac{1}{2}\left(\tau_k^\ast(X_k\otimes X_k) V_{k\ell}\tau_\ell
+\tau_\ell^\ast V_{\ell k}(X_k\otimes X_k)\tau_k\right)+\tau_k^\ast(X_\ell\otimes X_k)\tau_\ell.
\]
\noindent\textit{Proof. }  Under the above  assumptions,
\[
\Phi=\sum_{k=1}^K\tau_k^\ast V_k\tau_k=\sum_{k=1}^K\tau_k^\ast\tau_k=I_\mathcal{X},
\]
where $I_\mathcal{X}$ is the indentity operator of $\mathcal{X}$, and
\begin{eqnarray}\label{decompo}
\sqrt{n}\left(\widehat{T}_n-T\right)&=&\sqrt{n}\left(\widehat{\Phi}_n^{-1/2}\widehat{\Psi}_n\widehat{\Phi}_n^{-1/2}-\Psi\right)\nonumber\\
&=&\sqrt{n}\left(\widehat{\Phi}_n^{-1/2}-I_\mathcal{X}\right)\widehat{\Psi}_n\widehat{\Phi}_n^{-1/2}+\sqrt{n}\left(\widehat{\Psi}_n-\Psi\right)\widehat{\Phi}_n^{-1/2}+\Psi\left(\sqrt{n}(\widehat{\Phi}_n^{-1/2}-I_\mathcal{X})\right)\nonumber\\
&=&-\widehat{\Phi}_n^{-1}\left(\sqrt{n}(\widehat{\Phi}_n-I_\mathcal{X})\right)\left(\widehat{\Phi}_n^{-1/2}+I_\mathcal{X}\right)^{-1}\widehat{\Psi}_n\widehat{\Phi}_n^{-1/2}+\sqrt{n}\left(\widehat{\Psi}_n-\Psi\right)\widehat{\Phi}_n^{-1/2}\nonumber\\
&-&\Psi\widehat{\Phi}_n^{-1}\left(\sqrt{n}(\widehat{\Phi}_n-I_\mathcal{X})\right)\left(\widehat{\Phi}_n^{-1/2}+I_\mathcal{X}\right)^{-1}.
\end{eqnarray}
Clearly, 
\begin{equation}\label{rel1}
V_{k\ell}=\mathbb{E}(\tau_\ell(X)\otimes\tau_k( X))=\tau_kV\tau_\ell^\ast,
\end{equation}
where $V=\mathbb{E}(X\otimes X)$. Moreover, putting
\[
X^{(i)}=\left(
\begin{array}{c}
X_1^{(i)}\\
\vdots\\
X_K^{(i)}
\end{array}
\right),
\]
we have
\begin{eqnarray}\label{rel2}
\widehat{V}_{k\ell\cdot n}&=&\frac{1}{n}\sum_{i=1}^nX_\ell^{(i)}\otimes X_k^{(i)}-\overline{X}_{\ell\cdot n}\otimes\overline{X}_{k\cdot n}\nonumber\\
&=&\frac{1}{n}\sum_{i=1}^n\tau_\ell(X^{(i)})\otimes\tau_k( X^{(i)})-\tau_\ell(\overline{X}_ n)\otimes \tau_k(\overline{X}_n)\nonumber\\
&=&\tau_k\widehat{V}_n\tau_\ell^\ast,
\end{eqnarray}
where  $\overline{X}_ n=n^{-1}\sum_{i=1}^nX^{(i)}$ and
\begin{equation}\label{covemp}
\widehat{V}_n=\frac{1}{n}\sum_{i=1}^nX^{(i)}\otimes  X^{(i)}-\overline{X}_ n\otimes \overline{X}_n.
\end{equation}
Therefore, using  (\ref{rel1}) and (\ref{rel2}), we obtain
\begin{equation}\label{decompopsi}
\sqrt{n}\left(\widehat{\Psi}_n-\Psi\right)=\sum_{k=1}^K\sum_{\stackrel{\ell=1}{\ell\neq k}}^K\tau_k^\ast\tau_k \widehat{H}_n\tau_\ell^\ast\tau_\ell
=f(\widehat{H}_n),
\end{equation}
where $\widehat{H}_n=\sqrt{n}(\widehat{V}_{n}-V)$ and $f$  is the operator defined as
\[
f\,:\,A\in\mathcal{L}(\mathcal{X})\mapsto \sum_{k=1}^K\sum_{\stackrel{\ell=1}{\ell\neq k}}^K\tau_k^\ast\tau_k A\tau_\ell^\ast\tau_\ell\in\mathcal{L}(\mathcal{X}).
\]
Further, since $I_\mathcal{X}=\sum_{k=1}^K\tau_k^\ast\tau_k$, 
we obtain
\begin{equation}\label{decompophi}
\sqrt{n}(\widehat{\Phi}_n-I_\mathcal{X})=\sum_{k=1}^K\tau_k^\ast \tau_k\widehat{H}_{n}\tau_k^\ast\tau_k=g(\widehat{H}_{n}),
\end{equation}
where $g$ is the operator  $g\,:\,A\in\mathcal{L}(\mathcal{X})\mapsto \sum_{k=1}^K\tau_k^\ast\tau_k A\tau_k^\ast\tau_k\in\mathcal{L}(\mathcal{X})$. Then, using (\ref{decompo}), (\ref{decompopsi}) and (\ref{decompophi}),
we  obtain $\sqrt{n}\left(\widehat{T}_n-T\right)=\widehat{\varphi}_n(\widehat{H}_n)$, where $\widehat{\varphi}_n$ is the random operator from $\mathcal{L}(\mathcal{X})$ to itself defined by
\[
\widehat{\varphi}_n(A)= 
-(\widehat{\Phi}_n^{-1/2}+I_\mathcal{X})^{-1}g(A)\widehat{\Phi}_n^{-1}\widehat{\Psi}_n\widehat{\Phi}_n^{-1/2}+f(A)\widehat{\Phi}_n^{-1/2}-\Psi (\widehat{\Phi}_n^{-1/2}+I_\mathcal{X})^{-1}g(A)\widehat{\Phi}_n^{-1}.
\]
 Considering the operator
\[
\varphi\,:\,A\in\mathcal{L}(\mathcal{X})\mapsto
-\frac{1}{2}g(A)\Psi+f(A)-\frac{1}{2}\Psi g(A)\in\mathcal{L}(\mathcal{X}),
\]
and denoting by $\Vert\cdot\Vert_\infty$ (resp. $\Vert\cdot\Vert_{\infty\infty}$) the norm of $\mathcal{L}(\mathcal{X})$ (resp. $\mathcal{L}(\mathcal{L}(\mathcal{X}))$) defined by $\Vert A\Vert_\infty=\sup_{x\in \mathcal{X}-\{0\}}\Vert Ax\Vert_\mathcal{X}/\Vert x\Vert_\mathcal{X}$ (resp. $\Vert h\Vert_{\infty\infty}=\sup_{B\in \mathcal{L}(\mathcal{X})-\{0\}}\Vert h(B)\Vert_\infty/\Vert B\Vert_\infty$) for any $A$ (resp. $h$) in $\mathcal{L}(\mathcal{X})$ (resp. $\mathcal{L}(\mathcal{L}(\mathcal{X}))$), we have
\begin{eqnarray}\label{inegalites}
\Vert \widehat{\varphi}_n(\widehat{H}_{n})-\varphi(\widehat{H}_{n})\Vert_\infty
&= &\left\Vert -\left((\widehat{\Phi}_n^{-1/2}+I_\mathcal{X})^{-1}-\frac{1}{2}I_\mathcal{X}\right)g(\widehat{H}_{n})\widehat{\Phi}_n^{-1}\widehat{\Psi}_n\widehat{\Phi}_n^{-1/2}\right.\nonumber\\
&-&\frac{1}{2}g(\widehat{H}_{n})\left(\widehat{\Phi}_n^{-1}\widehat{\Psi}_n\widehat{\Phi}_n^{-1/2}-\Psi\right)
+ f(\widehat{H}_{n})(\widehat{\Phi}_n^{-1/2}-I_\mathcal{X})\nonumber\\
&-&\left.\Psi\left((\widehat{\Phi}_n^{-1/2}+I_\mathcal{X})^{-1}-\frac{1}{2}I_\mathcal{X}\right)g(\widehat{H}_{n})
\widehat{\Phi}_n^{-1}- \frac{1}{2}\Psi (\widehat{\Phi}_n^{-1}-I_\mathcal{X})\right\Vert_\infty\nonumber\\ 
&\leq &  \Vert (\widehat{\Phi}_n^{-1/2}+I_\mathcal{X})^{-1}-\frac{1}{2}I_\mathcal{X}\Vert_\infty\,\Vert g(\widehat{H}_{n})\Vert_\infty\,\Vert\widehat{\Phi}_n^{-1}\widehat{\Psi}_n\widehat{\Phi}_n^{-1/2}\Vert_\infty.\nonumber\\
&+&\frac{1}{2}\Vert g(\widehat{H}_{n})\Vert_\infty\,\Vert\widehat{\Phi}_n^{-1}\widehat{\Psi}_n\widehat{\Phi}_n^{-1/2}-\Psi\Vert_\infty
+ \Vert f(\widehat{H}_{n})\Vert_\infty\,\Vert\widehat{\Phi}_n^{-1/2}-I_\mathcal{X}\Vert_\infty\nonumber\\
&+&\Vert\Psi\Vert_\infty\,\Vert(\widehat{\Phi}_n^{-1/2}+I_\mathcal{X})^{-1}-\frac{1}{2}I_\mathcal{X}\Vert_\infty\,\Vert g(\widehat{H}_{n})
\Vert_\infty\,\Vert\widehat{\Phi}_n^{-1}\Vert_\infty\nonumber\\
& + &\frac{1}{2}\Vert\Psi\Vert_\infty\,\Vert g(\widehat{H}_{n})
\Vert_\infty\,\Vert \widehat{\Phi}_n^{-1}-I_\mathcal{X}\Vert_\infty\nonumber\\
&\leq &  \left(\Vert (\widehat{\Phi}_n^{-1/2}+I_\mathcal{X})^{-1}-\frac{1}{2}I_\mathcal{X}\Vert_\infty\,\Vert g\Vert_{\infty\infty}\,\Vert\widehat{\Phi}_n^{-1}\widehat{\Psi}_n\widehat{\Phi}_n^{-1/2}\Vert_\infty\right.\nonumber\\
&+&\frac{1}{2}\Vert g\Vert_{\infty\infty}\,\Vert\widehat{\Phi}_n^{-1}\widehat{\Psi}_n\widehat{\Phi}_n^{-1/2}-\Psi\Vert_\infty
+ \Vert f\Vert_{\infty\infty}\,\Vert\widehat{\Phi}_n^{-1/2}-I_\mathcal{X}\Vert_\infty\nonumber\\
&+&\Vert\Psi\Vert_\infty\,\Vert(\widehat{\Phi}_n^{-1/2}+I_\mathcal{X})^{-1}-\frac{1}{2}I_\mathcal{X}\Vert_\infty\,\Vert g\Vert_{\infty\infty}\,\Vert\widehat{\Phi}_n^{-1}\Vert_\infty\nonumber\\
& + &\left.\frac{1}{2}\Vert\Psi\Vert_\infty\,\Vert g\Vert_{\infty\infty}\,\Vert \widehat{\Phi}_n^{-1}-I_\mathcal{X}\Vert_\infty\right)\Vert\widehat{H}_{n}\Vert_\infty.
\end{eqnarray}
Using the strong law of large numbers, it is easy to verify that, for any $(k,\ell)\in\{1,\cdots,K\}^2$ with $k\neq \ell$,  $\widehat{V}_{k\ell\cdot n}$ (resp. $\widehat{V}_{k\cdot n}$) converge almost surely to $V_{k\ell}$ (resp. $\widehat{V}_{k}$), as $n\rightarrow +\infty$. Consequently, $\widehat{\Phi}_n$ (resp. $\widehat{\Psi}_n$) converge almost surely to $\Phi=I_\mathcal{X}$ (resp. $\Psi$), as $n\rightarrow +\infty$. This implies the almost sure convergence of $(\widehat{\Phi}_n^{-1/2}+I_\mathcal{X})^{-1}$ (resp.  $\widehat{\Phi}_n^{-1}\widehat{\Psi}_n\widehat{\Phi}_n^{-1/2}$; resp. $\widehat{\Phi}_n^{-1}$; resp. $\widehat{\Phi}_n^{-1/2}$) to $\frac{1}{2}I_\mathcal{X}$ (resp. $\Psi$; resp. $I_\mathcal{X}$; resp. $I_\mathcal{X}$), as $n\rightarrow +\infty$. Furthermore, denoting by $\Vert\cdot\Vert$ the norm of $\mathcal{L}(\mathcal{X})$ defined by $\Vert A\Vert=\sqrt{\textrm{tr}(AA^\ast)}$ and using the properties $(a\otimes b)(c\otimes d)=<a,d>c\otimes b$ and tr$(a\otimes b)=<a,b>$ of the tensor product (see Dauxois et al. (1994)), we have:
\begin{eqnarray*}
\mathbb{E}\left(\Vert X\otimes X\Vert^2\right)&=&\mathbb{E}\left(\textrm{tr}((X\otimes X)(X\otimes X)\right)=\mathbb{E}\left(\Vert  X\Vert^4_\mathcal{X}\right)
=\mathbb{E}\left(\left(\sum_{k=1}^K\Vert X_k\Vert_k^2\right)^2\right)\\
&=&\sum_{k=1}^K\mathbb{E}(\Vert X_k\Vert_k^4)+\sum_{k=1}^K\sum_{\stackrel{\ell=1}{\ell\neq k}}^K\mathbb{E}(\Vert X_k\Vert_k^2\Vert X_\ell\Vert_\ell^2)\\
&\leq &\sum_{k=1}^K\mathbb{E}(\Vert X_k\Vert_k^4)+\sum_{k=1}^K\sum_{\stackrel{\ell=1}{\ell\neq k}}^K\sqrt{\mathbb{E}(\Vert X_k\Vert_k^4)}\sqrt{\mathbb{E}(\Vert X_\ell\Vert_\ell^4)}<+\infty.
\end{eqnarray*}
Then,  the central limit theorem can be used. It gives the convergence in distribution, as $n\rightarrow +\infty$,  of $\sqrt{n}\left(n^{-1}\sum_{i=1}X^{(i)}\otimes X^{(i)}-V\right)$ to an random variable $H$ having the normal distribution in $\mathcal{L}(\mathcal{X})$ with mean equal to $0$ and a covariance operator equal to that of $X\otimes X$. Since, by the central limit theorem again, $\sqrt{n}\overline{X}_n$  converges in distribution, as $n\rightarrow +\infty$, to an random variable having a normal distribution in $\mathcal{X}$ with mean equal to $0$ and a covariance operator equal to  $V$, we deduce from the equality $\sqrt{n}\left(\overline{X}_n\otimes\overline{X}_n\right)=n^{-1/2}\left(\sqrt{n}\overline{X}_n\right)\otimes\left( \sqrt{n}\overline{X}_n\right)$ that  $\sqrt{n}\left(\overline{X}_n\otimes\overline{X}_n\right)$ converges in probability to $0$, as $n\rightarrow +\infty$. Therefore, from (\ref{covemp}) and Slutsky's theorem, we deduce that $\widehat{H}_n$ converges in distribution, as $n\rightarrow +\infty$ to  $H$. Then, from (\ref{inegalites}), we conclude that $\widehat{\varphi}_n(\widehat{H}_{n})-\varphi(\widehat{H}_{n})$ converges in probability to $0$, as $n\rightarrow +\infty$. Then, using again Slutsky's  theorem, we deduce that $\widehat{\varphi}_n(\widehat{H}_{n})$ and $\varphi(\widehat{H}_{n})$ both converge in distribution to the same distribution. Since $\varphi$ is a linear map (and is, therefore, continuous), this distribution   just is that of the random variable  $U=\varphi(H)$, that is the normal distribution in $\mathcal{L}(\mathcal{X})$ with mean $0$ and covariance operator equal to that of $Z=\varphi (X\otimes X)$. Clearly,
\[
g(X\otimes X)=\sum_{k=1}^K\tau_k^\ast\tau_k (X\otimes X)\tau_k^\ast\tau_k=\sum_{k=1}^K\tau_k^\ast((\tau_k (X))\otimes (\tau_k(X)))\tau_k=\sum_{k=1}^K\tau_k^\ast(X_k\otimes X_k)\tau_k,
\]
and
\[
f(X\otimes X)=\sum_{k=1}^K\sum_{\stackrel{\ell=1}{\ell\neq k}}^K\tau_k^\ast\tau_k (X\otimes X)\tau_\ell^\ast\tau_\ell
=\sum_{k=1}^K\sum_{\stackrel{\ell=1}{\ell\neq k}}^K\tau_k^\ast(X_\ell\otimes X_k)\tau_\ell.
\]
Then, since $\tau_k\tau_j^\ast=\delta_{kj}I_k$, it follows
\[
g(X\otimes X)\Psi=\sum_{k=1}^K\sum_{j=1}^K\sum_{\stackrel{\ell=1}{\ell\neq j}}^K\tau_k^\ast(X_k\otimes X_k)\tau_k\tau_j^\ast V_{j\ell}\tau_\ell
=\sum_{k=1}^K\sum_{\stackrel{\ell=1}{\ell\neq k}}^K\tau_k^\ast(X_k\otimes X_k) V_{k\ell}\tau_\ell
\]
and
\begin{eqnarray*}
\Psi g(X\otimes X)&=&\sum_{k=1}^K\sum_{\stackrel{\ell=1}{\ell\neq k}}^K\sum_{j=1}^K\tau_k^\ast V_{k\ell}\tau_\ell\tau_j^\ast(X_j\otimes X_j)\tau_j
=\sum_{k=1}^K\sum_{\stackrel{\ell=1}{\ell\neq k}}^K\tau_k^\ast V_{k\ell}(X_\ell\otimes X_\ell)\tau_\ell\\
&=&\sum_{k=1}^K\sum_{\stackrel{\ell=1}{\ell\neq k}}^K\tau_\ell^\ast V_{\ell k}(X_k\otimes X_k)\tau_k.
\end{eqnarray*}
Thus,
\[
Z
=\sum_{k=1}^K\sum_{\stackrel{\ell=1}{\ell\neq k}}^K-\frac{1}{2}\left(\tau_k^\ast(X_k\otimes X_k) V_{k\ell}\tau_\ell
+\tau_\ell^\ast V_{\ell k}(X_k\otimes X_k)\tau_k\right)+\tau_k^\ast(X_\ell\otimes X_k)\tau_\ell.
\]
\hfill$\Box$

\noindent Using the preceding theorem and results in Eaton and Tyler (1991,1994), we can now give asymptotic distributions for the canonical coefficients. We denote by $\left(\rho^\prime_j\right)_{1\leq j\leq r}$  (with $r\in\mathbb{N}^\ast$) the sequence of distinct eigienvalues of $T$ in decreasing order, that is  $\rho^\prime_1>\cdots>\rho^\prime_r$. Putting $m_0=0$, denoting by $m_j$ the multiplicity of $\rho^\prime_j$ and putting $\nu_j=\sum_{k=0}^{j-1}m_k$ for any  $j\in\{1,\cdots r\}$, it is clear that for any $i\in\{\nu_{j-1}+1,\cdots,\nu_j\}$ one has $\rho_i=\rho^\prime_j$.
Further, considering the eigenspace $E_j=\ker(T-\rho^\prime_jI)$, we have the following decomposition in orthogonal direct sum: $\mathcal{X}=E_1\oplus\cdots\oplus  E_r$. We denote by $\Pi_j$ the orthogonal projector from $\mathcal{X}$ onto $E_j$, and by  $\Delta $ the continuous map which associates to each self-adjoint operator $A$ the vector $\Delta (A)$ of its eigenvalues in nonincreasing order. For $j\in\{1,\cdots r\}$, we consider $m_j$-dimensional vector given by $\upsilon_j=\rho^\prime_j\mathbb{J}_{m_j}$, where $\mathbb{J}_q$ denotes the $q$-dimensional vector with elements all equal to $1$, and the $\mathbb{R}^{m_j}$- valued  random vector:
\[
\hat{\upsilon}^n_j=\left(
\begin{array}{c}
\widehat{\rho}_{\nu_{j-1}+1\cdot n}\\
\vdots\\
\widehat{\rho}_{\nu_{j}\cdot n}
\end{array}
\right).
\]
Then, putting
\[
\widehat{\Lambda}_n=\left(
\begin{array}{c}
\hat{\upsilon}^n_1\\
\vdots\\
\hat{\upsilon}^n_r
\end{array}
\right)\,\,\,\,\textrm{and}\,\,\,\,
\Lambda=\left(
\begin{array}{c}
\upsilon_1\\
\vdots\\
\upsilon_r
\end{array}
\right),
\]
we have:

\bigskip

\noindent\textbf{Theorem 3.3. }\textsl{  $\sqrt{n}\left(\widehat{\Lambda}_n -\Lambda\right)$ converges in distribution, as $n\rightarrow +\infty$, to the $\mathbb{R}^p$-valued random vector}
\begin{eqnarray}\label{loivalp}
\zeta=\left(
\begin{array}{c}
\Delta (\Pi_1W\Pi_1)\\
\vdots\\
\Delta (\Pi_rW\Pi_r)
\end{array}
\right),
\end{eqnarray}
\textsl{where $W$ is  a random variable having a normal distribution in $\mathcal{L}(\mathcal{X})$, with mean $0$ and covariance operator $\Theta$ given by:}
\[
\Theta=\sum_{1\leq m,r,s,t\leq p}C(m,r,s,t)\,\,(e_m\otimes e_r)\widetilde{\otimes}(e_s\otimes e_t)
\]
\textsl{with}
\begin{eqnarray*}
C(m,r,s,t)=\sum_{k=1}^K\sum_{j=1}^K\sum_{\stackrel{\ell=1}{\ell\neq k}}^K\sum_{\stackrel{q=1}{q\neq j}}^K & &\left(\gamma^{m,r,s,t}_{k\ell jq}+\gamma^{m,r,t,s}_{k\ell jq}+\gamma^{r,m,s,t}_{k\ell jq}+\gamma^{r,m,t,s}_{k\ell jq}\right.\\
 & &\left. -\theta^{m,r,s,t}_{k\ell jq}-\theta^{r,m,s,t}_{k\ell jq}-\theta^{s,t,m,r}_{k\ell jq}-\theta^{t,s,m,r}_{k\ell jq}+\lambda^{m,r,s,t}_{k\ell jq}\right),
\end{eqnarray*}
\[
\gamma^{a,b,c,d}_{k\ell jq}=\frac{1}{4}\mathbb{E}\left(<X_k,\tau_k\beta^{(a)}>_k<X_k,V_{k\ell}\tau_\ell\beta^{(b)}>_k<X_j,\tau_j\beta^{(c)}>_j<X_j,V_{jq}\tau_q\beta^{(d)}>_j\right),
\]
\[
\theta^{a,b,c,d}_{k\ell jq}=\frac{1}{2}\mathbb{E}\left(<X_k,\tau_k\beta^{(a)}>_k<X_k,V_{k\ell}\tau_\ell\beta^{(b)}>_k<X_j,\tau_j\beta^{(c)}>_j<X_q,\tau_q\beta^{(d)}>_q\right)
\]
\textsl{and}
\[
\gamma^{a,b,c,d}_{k\ell jq}=\mathbb{E}\left(<X_k,\tau_k\beta^{(a)}>_k<X_\ell,\tau_\ell\beta^{(b)}>_\ell<X_j,\tau_j\beta^{(c)}>_j<X_q,\tau_q\beta^{(d)}>_q\right).
\]
\noindent\textit{Proof. } Since $\Delta (\widehat{T}_n)=\widehat{\Lambda}_n$ and  $\Delta (T)=\Lambda$, we deduce from Theorem 3.2 and the Theorem 2.1 of Eaton and Tyler (1994)  that  $\sqrt{n}\left(\widehat{\Lambda}_n -\Lambda\right)$ converges in distribution, as $n\rightarrow +\infty$, to the random variable given in (\ref{loivalp}) with  $W=P^\ast UP$, where  $P=\sum_{\ell=1}^pe_\ell\otimes \beta^{(\ell )}$. Clearly, $W$ has a normal distribution with mean $0$ and covariance operator $\Theta$ equal to that of $P^\ast ZP$. In order to give an explicit expression of $\Theta$, let us first note that:
\begin{eqnarray*}
P^\ast ZP& =&\sum_{k=1}^K\sum_{\stackrel{\ell=1}{\ell\neq k}}^K-\frac{1}{2}\left(P^\ast\tau_k^\ast(X_k\otimes X_k) V_{k\ell}\tau_\ell P
+P^\ast\tau_\ell^\ast V_{\ell k}(X_k\otimes X_k)\tau_kP\right)\\
& &+P^\ast\tau_k^\ast(X_\ell\otimes X_k)\tau_\ell P\\
&=&\sum_{k=1}^K\sum_{\stackrel{\ell=1}{\ell\neq k}}^K-\frac{1}{2}\left((P^\ast \tau_\ell^\ast  V_{\ell k}X_k)\otimes(P^\ast\tau_k^\ast X_k)
+(P^\ast\tau_k^\ast X_k)\otimes ( P^\ast\tau_\ell^\ast V_{\ell k}X_k)\right)\\
& &+(P^\ast\tau_\ell^\ast X_\ell)\otimes (P^\ast\tau_k^\ast X_k).
\end{eqnarray*}
Since
\begin{eqnarray*}
P^\ast \tau_\ell^\ast  V_{\ell k}X_k&=&\left(\sum_{m=1}^p \beta^{(m )}\otimes e_m\right)\tau_\ell^\ast  V_{\ell k}X_k
=\sum_{m=1}^p <\beta^{(m)},\tau_\ell^\ast  V_{\ell k}X_k>_\mathcal{X}e_m\\
&=&\sum_{m=1}^p <\tau_\ell\beta^{(m)},V_{\ell k}X_k>_\ell e_m
\end{eqnarray*}
and, similarly, $
P^\ast \tau_k^\ast X_k=\sum_{m=1}^p <\tau_k\beta^{(m)},X_k>_ke_m$, it follows:
\begin{eqnarray*}
P^\ast ZP=\sum_{m=1}^p\sum_{r=1}^p
\left[\sum_{k=1}^K\sum_{\stackrel{\ell=1}{\ell\neq k}}^K \right.&-& \frac{1}{2}(<\tau_\ell\beta^{(m)},V_{\ell k}X_k>_\ell<\tau_k\beta^{(r)},X_k>_k\\
&+&<\tau_\ell\beta^{(r)},V_{\ell k}X_k>_\ell<\tau_k\beta^{(m)},X_k>_k)\\
& +&\left.<\tau_\ell\beta^{(m)},X_\ell>_\ell<\tau_k\beta^{(r)},X_k>_k\right]\,e_m\otimes e_r.
\end{eqnarray*}
From:
\begin{eqnarray*}
\mathbb{E}\left(<\tau_\ell\beta^{(m)},V_{\ell k}X_k>_\ell<\tau_k\beta^{(r)},X_k>_k\right)&=&\mathbb{E}\left(<(X_k\otimes X_k)(\tau_k\beta^{(r)}),V_{k\ell}\tau_\ell\beta^{(m)}>_k\right)\\
&=&<\mathbb{E}(X_k\otimes X_k)(\tau_k\beta^{(r)}),V_{k\ell}\tau_\ell\beta^{(m)}>_k\\
&=&<V_k\tau_k\beta^{(r)},V_{k\ell}\tau_\ell\beta^{(m)}>_k\\
&=&<\tau_k\beta^{(r)},V_{k\ell}\tau_\ell\beta^{(m)}>_k,
\end{eqnarray*}
\begin{eqnarray*}
\mathbb{E}\left(<\tau_\ell\beta^{(r)},V_{\ell k}X_k>_\ell<\tau_k\beta^{(m)},X_k>_k\right)=<\tau_k\beta^{(m)},V_{k\ell}\tau_\ell\beta^{(r)}>_k
\end{eqnarray*}
and
\begin{eqnarray*}
\mathbb{E}\left(<\tau_\ell\beta^{(m)},X_\ell>_\ell<\tau_k\beta^{(r)},X_k>_k\right)&=&\mathbb{E}\left(<(X_\ell\otimes X_k)(\tau_\ell\beta^{(m)}),\tau_k\beta^{(r)}>_k\right)\\
&=&<\mathbb{E}(X_\ell\otimes X_k)(\tau_\ell\beta^{(m)}),\tau_k\beta^{(r)}>_k\\
&=&<V_{k\ell}\tau_\ell\beta^{(m)},\tau_k\beta^{(r)}>_k,
\end{eqnarray*}
we deduce that $\mathbb{E}(P^\ast ZP)=0$. Thus,
\begin{eqnarray*}
\Theta&=&\mathbb{E}\left((P^\ast ZP)\widetilde{\otimes}(P^\ast ZP)\right)=\sum_{1\leq m,r,s,t\leq p}C(m,r,s,t)\,\,(e_m\otimes e_r)\widetilde{\otimes}(e_s\otimes e_t),
\end{eqnarray*}
where
\begin{eqnarray*}
C(m,r,s,t)=\sum_{k=1}^K\sum_{j=1}^K\sum_{\stackrel{\ell=1}{\ell\neq k}}^K\sum_{\stackrel{q=1}{q\neq j}}^K &  & \mathbb{E}\left(Y^{m,r}_{k\ell}Y^{s,q}_{jq}\right)
\end{eqnarray*}
with
\begin{eqnarray*}
Y^{m,r}_{k\ell}= & &-\frac{1}{2}\left(<\tau_\ell\beta^{(m)},V_{\ell k}X_k>_\ell<\tau_k\beta^{(r)},X_k>_k
+<\tau_\ell\beta^{(r)},V_{\ell k}X_k>_\ell<\tau_k\beta^{(m)},X_k>_k\right)\\
 & &+<\tau_\ell\beta^{(m)},X_\ell>_\ell<\tau_k\beta^{(r)},X_k>_k.
\end{eqnarray*}
Further calculations give
 \begin{eqnarray*}
 \mathbb{E}\left(Y^{m,r}_{k\ell}Y^{s,q}_{jq}\right)=& &\gamma^{m,r,s,t}_{k\ell jq}+\gamma^{m,r,t,s}_{k\ell jq}+\gamma^{r,m,s,t}_{k\ell jq}+\gamma^{r,m,t,s}_{k\ell jq}\\
 & &  -\theta^{m,r,s,t}_{k\ell jq}-\theta^{r,m,s,t}_{k\ell jq}-\theta^{s,t,m,r}_{k\ell jq}-\theta^{t,s,m,r}_{k\ell jq}+\lambda^{m,r,s,t}_{k\ell jq}.
\end{eqnarray*}
\hfill $\Box$

When $T$ has simple eigenvalues, that is $\rho_1>\rho_2>\cdots >\rho_q$, the  preceding theorem has a simpler statement. We have:

\bigskip

\noindent\textbf{Corollary 3.1. }\textsl{ When the eigenvalues of  $T$ are simple,  $\sqrt{n}\left(\widehat{\Lambda}_n -\Lambda\right)$ converges in distribution, as $n\rightarrow +\infty$, to a random variable having a normal distribution in  $\mathbb{R}^p$ with mean $0$ and covariance matrix $\Sigma=\left(\sigma_{ij}\right)_{1\leq i,j\leq p}$ with:}
\[
\sigma_{ij}=\sum_{1\leq m,r,s,t\leq p}\beta^{(i)}_m\beta^{(i)}_r\beta^{(j)}_s\beta^{(j)}_tC(m,r,s,t).
\]

\noindent\textit{Proof. }In this case, $m_1=\cdots m_p=1$ and, for any $j\in\{1,\cdots,p\}$, $\Pi_j=\beta^{(j)}\otimes\beta^{(j)}$. Thus
\begin{eqnarray*}
\Pi_jW\Pi_j&=&((\beta^{(j)}\otimes\beta^{(j)})W(\beta^{(j)}\otimes\beta^{(j)})
=(\beta^{(j)}\otimes\beta^{(j)})(\beta^{(j)}\otimes(W\beta^{(j)}))\\
&=&<\beta^{(j)},W\beta^{(j)}>_\mathcal{X}\beta^{(j)}\otimes\beta^{(j)},
\end{eqnarray*}
and, therefore, $\Delta (\Pi_jW\Pi_j)=<\beta^{(j)},W\beta^{(j)}>_\mathcal{X}$. Then, $\zeta$ is a linear function of $W$ and, consequently, it has a normal distribution with mean $0$ and covariance matrix $\Sigma=\left(\sigma_{ij}\right)_{1\leq i,j\leq p}$ with $\sigma_{ij}=\mathbb{E}\left(<\beta^{(i)},W\beta^{(i)}>_\mathcal{X}<\beta^{(j)},W\beta^{(j)}>_\mathcal{X}\right)$. Denoting by $<\cdot , \cdot >$ the inner product of operators defined by $<A,B>=\textrm{tr }(AB^\ast)$, we have:
\[
<W,\beta^{(j)}\otimes\beta^{(j)}>=\textrm{tr}\left(W(\beta^{(j)}\otimes\beta^{(j)})\right)=\textrm{tr}\left(\beta^{(j)}\otimes (W\beta^{(j)})\right)
=<\beta^{(j)},W\beta^{(j)}>_\mathcal{X},
\]
it follows that
\begin{eqnarray*}
\sigma_{ij}&=&\mathbb{E}\left(<\beta^{(i)},W\beta^{(i)}>_\mathcal{X}<\beta^{(j)},W\beta^{(j)}>_\mathcal{X}\right)\\
&=&\mathbb{E}\left(<W,\beta^{(i)}\otimes\beta^{(i)}><W,\beta^{(j)}\otimes\beta^{(j)}>\right)\\
&=&\mathbb{E}\left(<(W\widetilde{\otimes}W)(\beta^{(i)}\otimes\beta^{(i)}),\beta^{(j)}\otimes\beta^{(j)}>\right)\\
&=&<\mathbb{E}(W\widetilde{\otimes}W)(\beta^{(i)}\otimes\beta^{(i)}),\beta^{(j)}\otimes\beta^{(j)}>\\
&=&<\Theta (\beta^{(i)}\otimes\beta^{(i)}),\beta^{(j)}\otimes\beta^{(j)}>\\
&=&\sum_{1\leq m,r,s,t\leq p}C(m,r,s,t)<\left((e_m\otimes e_r)\widetilde{\otimes}(e_s\otimes e_t)\right) (\beta^{(i)}\otimes\beta^{(i)}),\beta^{(j)}\otimes\beta^{(j)}>\\
&=&\sum_{1\leq m,r,s,t\leq p}C(m,r,s,t)<e_m\otimes e_r,\beta^{(i)}\otimes\beta^{(i)}><e_s\otimes e_t,\beta^{(j)}\otimes\beta^{(j)}>.
\end{eqnarray*}
Then, the  required result is obtained from
\begin{eqnarray*}
<e_m\otimes e_r,\beta^{(i)}\otimes\beta^{(i)}>&=&\textrm{tr}\left((e_m\otimes e_r)(\beta^{(i)}\otimes\beta^{(i)}>)\right)\\
&=&\textrm{tr}\left(<e_m,\beta^{(i)}>_\mathcal{X}\beta^{(i)}\otimes e_r\right)\\
&=&<e_m,\beta^{(i)}>_\mathcal{X}<e_r,\beta^{(i)}_\mathcal{X}>\\
&=&\beta^{(i)}_m\beta^{(i)}_r
\end{eqnarray*}
and $<e_s\otimes e_t,\beta^{(j)}\otimes\beta^{(j)}>=\beta^{(j)}_s\beta^{(j)}_t$.
\hfill $\Box$

\section{Testing for mutual non-correlation}
In this section, we consider the problem of testing for mutual non-correlation between $X_1,X_2,\cdots,X_K$, that is   testing for non-correlation between any pair $(X_k,X_\ell)$ for $ (k,\ell)\in\{1,\cdots,K\}^2$. The null  hypothesis  is
\[
\mathcal{H}_0\,:\,\forall (k,\ell)\in\{1,\cdots,K\}^2,\,k\neq\ell,\,\,V_{k\ell}=0
\]
and the alternative is given by:
\[
\mathcal{H}_1\,:\,\exists (k,\ell)\in\{1,\cdots,K\}^2,\,k\neq\ell,\,\,V_{k\ell}\neq 0,
\]
that is  testing for non-correlation between any pair $(X_k,X_\ell)$ of variables. We will first introduce a  test statistic for this problem, then  we will derive its asymptotic distribution under the null hypothesis, in the general case and in case $X$ has an elliptical distribution, by using the results of the preceding section.
\subsection{A test statistic}
For $(k,\ell)\in\{1,\cdots,K\}^2$,  denoting by $\pi_{k\ell}$ the operator 
\[
\pi_{k\ell}\,:\,A\in\mathcal{L}(\mathcal{X})\mapsto \tau_k A\tau_\ell^\ast\in\mathcal{L}(\mathcal{X}_\ell,\mathcal{X}_k),
\] 
we have for  $k\neq\ell$:
\begin{eqnarray*}
\pi_{k\ell}(T)&=&\pi_{k\ell}(\Phi^{-1/2}\Psi\Phi^{-1/2})=\sum_{i=1}^K\sum_{\stackrel{j=1}{j\neq i}}^K\tau_k\tau_i^\ast V_i^{-1/2} V_{ij}V_{j}^{-1/2}\tau_j\tau_\ell^\ast\\
&=&\sum_{i=1}^K\sum_{\stackrel{j=1}{j\neq i}}^K\delta_{ki}\delta_{j\ell} V_i^{-1/2} V_{ij}V_{j}^{-1/2}
=V_k^{-1/2} V_{k\ell}V_{\ell}^{-1/2}.
\end{eqnarray*}
Therefore, $\mathcal{H}_0$ is equivalent to having
$
\sum_{k=2}^K\sum_{\ell=1}^{k-1}\textrm{tr}\left(\pi_{k\ell}(T)\pi_{k\ell}(T)^\ast\right)=0$. This leads us to take as test statistic the random variable  $\widehat{S}_n$ given by:
\[
\widehat{S}_n=\sum_{k=2}^K\sum_{\ell=1}^{k-1}\textrm{tr}\left(\pi_{k\ell}(\widehat{T}_n)\pi_{k\ell}(\widehat{T}_n)^\ast\right).
\]
\subsection{Asymptotic distribution under null hypothesis}
 For $k\in\{1,\cdots, K\}$, we denote by $\{e_j^{(k)};\,j=1,\cdots, p_k\}$ an orthonormal basis of $\mathcal{X}_k$, where $p_k=\dim(\mathcal{X}_k)$, and we consider the matrix:
\begin{equation}\label{cov}
\Gamma=\left(
\begin{array}
[c]{ccccc}%
\Gamma_{21,21} & \Gamma_{21,31} & \Gamma_{21,32}  &\cdots & \Gamma_{21,KK-1}\\
\Gamma_{31,21} & \Gamma_{31,31} & \Gamma_{31,32}&\cdots & \Gamma_{31,KK-1}\\
\vdots & \vdots & \vdots & \cdots&\vdots\\
\Gamma_{KK-1,21} & \Gamma_{KK-1,31} & \Gamma_{KK-1,32}&\cdots & \Gamma_{KK-1,KK-1}
\end{array}
\right) ,
\end{equation}
where $\Gamma_{k\ell,rs}$ is the $p_kp_\ell\times p_rp_s$ matrix given by
\begin{equation}\label{gamma}
\Gamma_{k\ell,rs}=\left(
\begin{array}
[c]{ccccccccc}%
\gamma_{1111}^{k\ell,rs} & \gamma_{1121}^{k\ell,rs} & \cdots & \gamma_{11p_r1}^{k\ell,rs}&\cdots &\gamma_{111p_s}^{k\ell,rs} & \gamma_{112p_s}^{k\ell,rs} & \cdots & \gamma_{11p_rp_s}^{k\ell,rs}\\
\gamma_{2111}^{k\ell,rs} & \gamma_{2121}^{k\ell,rs} & \cdots & \gamma_{21p_r1}^{k\ell,rs}&\cdots &\gamma_{211p_s}^{k\ell,rs} & \gamma_{212p_s}^{k\ell,rs} & \cdots & \gamma_{21p_rp_s}^{k\ell,rs}\\
\vdots &\vdots &\cdots &\vdots &\cdots &\vdots & \vdots &\cdots &\vdots \\
\gamma_{p_k111}^{k\ell,rs} & \gamma_{p_k121}^{k\ell,rs} & \cdots & \gamma_{p_k1p_r1}^{k\ell,rs}&\cdots &\gamma_{p_k11p_s}^{k\ell,rs} & \gamma_{p_k12p_s}^{k\ell,rs} & \cdots & \gamma_{p_k1p_rp_s}^{k\ell,rs}\\
\vdots &\vdots &\cdots &\vdots &\cdots &\vdots & \vdots &\cdots &\vdots \\
\gamma_{1p_\ell 11}^{k\ell ,rs} & \gamma_{1p_\ell 21}^{k\ell ,rs} & \cdots & \gamma_{1p_\ell p_r1}^{k\ell ,rs}&\cdots &\gamma_{1p_\ell 1p_s}^{k\ell ,rs} & \gamma_{1p_\ell 2p_s}^{k\ell ,rs} & \cdots & \gamma_{1p_\ell p_rp_s}^{k\ell ,rs}\\
\gamma_{2p_\ell 11}^{k\ell,rs} & \gamma_{2p_\ell 21}^{k\ell ,rs} & \cdots & \gamma_{2p_\ell p_r1}^{k\ell,rs}&\cdots &\gamma_{2p_\ell 1p_s}^{k\ell,rs} & \gamma_{2p_\ell 2p_s}^{k\ell,rs} & \cdots & \gamma_{2p_\ell p_rp_s}^{k\ell,rs}\\
\vdots &\vdots &\cdots &\vdots &\cdots &\vdots & \vdots &\cdots &\vdots \\
\gamma_{p_kp_\ell 11}^{k\ell,rs} & \gamma_{p_kp_\ell 21}^{k\ell,rs} & \cdots & \gamma_{p_kp_\ell p_r1}^{k\ell,rs}&\cdots &\gamma_{p_kp_\ell 1p_s}^{k\ell,rs} & \gamma_{p_kp_\ell 2p_s}^{k\ell,rs} & \cdots & \gamma_{p_kp_\ell p_rp_s}^{k\ell,rs}
\end{array}
\right) ,
\end{equation}
with
\begin{eqnarray*}
\gamma^{k\ell,rs}_{ijpq}= \mathbb{E}\left( <X_k,e_i^{(k)}>_k\,<X_r,e_p^{(r)}>_r\,<X_\ell,e_j^{(\ell)}>_\ell\,<X_s,e_q^{(s)}>_s\right).
\end{eqnarray*}
Then, putting $d=\sum_{k=1}^K\sum_{\ell =1}^{k-1}p_kp_\ell$, we have:

\bigskip

\noindent\textbf{Theorem 4.1.}\textsl{ Under $\mathcal{H}_0$, the sequence $n\widehat{S}_n$ converges in distribution, as $n\rightarrow +\infty$,  to $\mathcal{Q}=\mathbb{W}^T\mathbb{W}$, where $\mathbb{W}$ is a random variable having a centered normal distribution in $\mathbb{R}^d$  with covariance matrix $\Gamma$.}

\noindent\textit{Proof. } First, note that under $\mathcal{H}_0$ we have $\Psi=0$ and, therefore, $T=0$. Then, the asymptotic distribution of $\widehat{S}_n$ under $\mathcal{H}_0$ can be obtained Theorem 3.2 which shows  that $\sqrt{n}\widehat{T}_n$ converges in distribution, as $n\rightarrow +\infty$, to a r.v. $U$ having a normal distribution in $\mathcal{L}(\mathcal{X})$ with mean $0$ and covariance operator equal to that of $Z$. Then, since the map $A\in\mathcal{L}(\mathcal{X})\mapsto\sum_{k=2}^K\sum_{\ell=1}^{k-1}\textrm{tr}\left(\pi_{k\ell}(A)\pi_{k\ell}(A)^\ast\right)\in\mathbb{R}$ is continuous, we deduce that $n\widehat{S}_n$  converges in distribution, as $n\rightarrow +\infty$, to $\mathcal{Q}=\sum_{k=2}^K\sum_{\ell=1}^{k-1}\textrm{tr}\left(W_{k\ell}W_{k\ell}^\ast\right)$, where $W_{k\ell}=\pi_{k\ell}(U)$. We can write $W_{k\ell}=\sum_{i=1}^{p_k}\sum_{j=1}^{p_\ell}W_{k\ell}^{ij}\,e_j^{(\ell)}\otimes e_i^{(k)}$ where  $W_{k\ell}^{ij}=<W_{k\ell},e_j^{(\ell)}\otimes e_i^{(k)}>$, and, using the equalities $(u\otimes v)(w\otimes z)=<u,z>w\otimes v$ and $\textrm{tr}(u\otimes v)=<u,v>$ (cf. Dauxois et al. (1994)), we have: 
\[
\textrm{tr}(W_{k\ell}W_{k\ell}^\ast)=\sum_{i=1}^{p_k}\sum_{j=1}^{p_\ell}\sum_{r=1}^{p_k}\sum_{s=1}^{p_\ell}\delta_{js}\delta_{ir}W_{k\ell}^{ij}W_{k\ell}^{rs}
=\sum_{i=1}^{p_k}\sum_{j=1}^{p_\ell}\left(W_{k\ell}^{ij}\right)^2
=\mathbb{W}_{k,\ell}^T\mathbb{W}_{k,\ell}.
\]
where $\delta$ denotes the usual Kronecker symbol and 
\[
\mathbb{W}_{k,\ell}=\left(W_{k\ell}^{11},W_{k\ell}^{21},\cdots,W_{k\ell}^{p_k1},\cdots,W_{k\ell}^{1p_\ell},W_{k\ell}^{2p_\ell},\cdots,W_{k\ell}^{p_kp_\ell}\right)^T.
\] 
Hence, $\mathcal{Q}=\mathbb{W}^T\mathbb{W}$, where $\mathbb{W}=\left(\mathbb{W}_{2,1}^T,\mathbb{W}_{3,1}^T,\mathbb{W}_{3,2}^T,\cdots,\mathbb{W}_{K,1}^T,\cdots,\mathbb{W}_{K,K-1}^T\right)^T$. Since $\mathbb{W}$ is a linear function of $U$, it is also normally distributed with mean $0$ and covariance matrix $\Gamma $ defined in  (\ref{cov}) and (\ref{gamma}) with $\gamma^{k\ell,rs}_{ijpq}=\mathbb{E}\left(W_{k\ell}^{ij}W_{rs}^{pq}\right)$. Further,
\begin{eqnarray*}
\gamma^{k\ell,rs}_{ijpq}&=& \mathbb{E}\left(<W_{k\ell},e_j^{(\ell)}\otimes e_i^{(k)}><W_{rs},e_q^{(s)}\otimes e_p^{(r)}>\right)\\
&=& \mathbb{E}\left(<(W_{k\ell}\widetilde{\otimes}W_{rs})(e_j^{(\ell)}\otimes e_i^{(k)}),e_q^{(s)}\otimes e_p^{(r)}>\right)\\
&=& <\mathbb{E}(W_{k\ell}\widetilde{\otimes}W_{rs})(e_j^{(\ell)}\otimes e_i^{(k)}),e_q^{(s)}\otimes e_p^{(r)}>,
\end{eqnarray*}
where $\widetilde{\otimes}$ denotes the tensor product related to the inner product of operators $<A,B>=$tr $(AB^\ast)$. Recalling that  $U$ has the same covariance operator than $\varphi(X\otimes X)$, that is  $\mathbb{E}(U\widetilde{\otimes}U)=\mathbb{E}\left((\varphi(X\otimes X))\widetilde{\otimes}(\varphi(X\otimes X))\right)$, we obtain 
\begin{eqnarray*}
\mathbb{E}\left(W_{k\ell}\widetilde{\otimes}W_{rs}\right)&=&\mathbb{E}\left((\pi_{k\ell}(U))\widetilde{\otimes}(\pi_{rs}(U))\right)
=\pi_{rs}\mathbb{E}\left(U\widetilde{\otimes}U\right)\pi_{k\ell}^\ast\\
&=&\pi_{rs}\mathbb{E}\left((\varphi(X\otimes X))\widetilde{\otimes}(\varphi(X\otimes X))\right)\pi_{k\ell}^\ast\\
&=&\mathbb{E}\left((\pi_{k\ell}(\varphi(X\otimes X)))\widetilde{\otimes}(\pi_{rs}(\varphi(X\otimes X)))\right)\\
&=&\mathbb{E}\left((X_\ell\otimes X_k)\widetilde{\otimes}(X_s\otimes X_r)\right).
\end{eqnarray*}
Thus, $\mathbb{E}\left(\mathbb{W}_{k\ell}\mathbb{W}_{rq}^T\right)=$, $\mathbb{E}\left(<X_\ell\otimes X_k,e_j^{(\ell)}\otimes e_i^{(k)}><X_s\otimes X_r,e_q^{(s)}\otimes e_p^{(r)}>\right)$
\begin{eqnarray*}
\gamma^{k\ell,rs}_{ijpq}
&=& \mathbb{E}\left(<X_\ell\otimes X_k,e_j^{(\ell)}\otimes e_i^{(k)}><X_s\otimes X_r,e_q^{(s)}\otimes e_p^{(r)}>\right)\\
&=& \mathbb{E}\left( <X_k,e_i^{(k)}>_k\,<X_r,e_p^{(r)}>_r\,<X_\ell,e_j^{(\ell)}>_\ell\,<X_s,e_q^{(s)}>_s\right).
\end{eqnarray*}
\hfill $\Box$
\subsection{Case of elliptical distribution}

\noindent Recall that a random variable valued into an Euclidean space $E$ is said to have an elliptical distribution $\mathcal{E}(\mu,\Sigma, h)$ if its characteristic function is of the form: 
\[
\exp(i<t,\mu>_E)\,h\left(<t,\Sigma t>_E\right).
\] 
In this subsection, we assume that  $X$ has an elliptical distribution  $\mathcal{E}(0,V,h)$. Then, the previous theorem has a more explicit formulation stated below.

\bigskip

\noindent\textbf{Theorem 4.2.}\textsl{When $X$ has the elliptical distribution $\mathcal{E}(0,V,h)$  then,  under $\mathcal{H}_0$, the sequence $n\widehat{S}_n$ converges in distribution, as $n\rightarrow +\infty$,  to $4h^{\prime\prime}(0)\chi^2_d$.}
\noindent\textit{Proof. } The $\mathbb{R}^4$- valued  random vector
\[
S=\left(
\begin{array}{c}
<X_k,e_i^{(k)}>_k\\
<X_r,e_p^{(r)}>_r\\
<X_\ell,e_j^{(\ell)}>_\ell\\
<X_s,e_q^{(s)}>_s
\end{array}
\right)
\]
can obviously been writen as $S=L(X)$, where $L$ in a linear map from $\mathcal{X}$ to $\mathbb{R}^4$. Then, $S$ has the elliptical distribution $E(0,LVL^\ast,h)$ and, therefore,
\begin{eqnarray*}
\gamma^{k\ell,rs}_{ijpq}&=&4h^{\prime\prime}(0)\left[\mathbb{E}\left( <X_k,e_i^{(k)}>_k\,<X_r,e_p^{(r)}>_r\right)\mathbb{E}\left( \,<X_\ell,e_j^{(\ell)}>_\ell\,<X_s,e_q^{(s)}>_s\right)\right.\\
& &\,\,\,\,\,\,\,\,\,\,\,\,\,\,\,\,\,\,+\mathbb{E}\left( <X_k,e_i^{(k)}>_k\,<X_\ell,e_j^{(\ell)}>_\ell\right)\mathbb{E}\left( <X_r,e_p^{(r)}>_r\,<X_s,e_q^{(s)}>_s\right)\\
& & \,\,\,\,\,\,\,\,\,\,\,\,\,\,\,\,\,\,\left.+\mathbb{E}\left( <X_k,e_i^{(k)}>_k\,<X_s,e_q^{(s)}>_s\right)\mathbb{E}\left( <X_r,e_p^{(r)}>_r\,<X_\ell,e_j^{(\ell)}>_\ell\right)\right]\\
&=&4h^{\prime\prime}(0)\left[<V_{kr}e^{(r)}_p,e^{(k)}_i>_k<V_{\ell s}e^{(s)}_q,e^{(\ell)}_j>_\ell \right]\\
& & \,\,\,\,\,\,\,\,\,\,\,\,\,\,\,\,\,\,\,\,\,+ <V_{k\ell}e^{(\ell)}_j,e^{(k)}_i>_k<V_{r s}e^{(s)}_q,e^{(r)}_p>_p\\
& & \,\,\,\,\,\,\,\,\,\,\,\,\,\,\,\,\,\,\,\,\,+ \left.<V_{k s}e^{(s)}_q,e^{(k)}_i>_k<V_{\ell r}e^{(r)}_p,e^{(\ell)}_j>_\ell\right].
\end{eqnarray*}
Since, under $\mathcal{H}_0$, we have for any $(m,l)\in\{1,\cdots,K\}^2$, $V_{m l}=\delta_{m l}V_l=\delta_{m l}I_l$, where $\delta$ denotes the usual Kronecker symbol, we finally obtain:
\begin{eqnarray}\label{gamfinal}
\gamma^{k\ell,rs}_{ijpq}=4h^{\prime\prime}(0)\left(\delta_{kr}\delta_{\ell s}\delta_{ip}\delta_{jq}+\delta_{k\ell}\delta_{r s}\delta_{ij}\delta_{pq}+\delta_{ks}\delta_{\ell r}\delta_{iq}\delta_{jp}\right).
\end{eqnarray}
If $k=r$, $\ell =s$ with $k\in\{2,\cdots, K\}$ and $\ell\in\{1,\cdots,k-1\}$, since $\ell\neq k$ and, equivalently, $s\neq r$, it follows that  
\[
\gamma^{k\ell,rs}_{ijpq}=4h^{\prime\prime}(0)\delta_{ip}\delta_{jq}=\left\{
\begin{array}{ccl}
4h^{\prime\prime}(0) & &\textrm{if }i=j=p=q \\
 & & \\
0 & & \textrm{otherwise}
\end{array}
\right. .
\]
Hence, 
\begin{eqnarray}\label{gamdiag}
\Gamma_{k\ell,k\ell}=4h^{\prime\prime}(0)\,\mathbb{I}_{p_kp_\ell},
\end{eqnarray} 
where $\mathbb{I}_m$ denotes the $m$-dimensional identity matrix.  If $(k,\ell)\neq (r,s)$ with $(k,r)\in\{2,\cdots,K\}^2$ and $(\ell,s)\in\{1,\cdots,k-1\}\times \{1,\cdots ,r -1\}$, suppose that  $\gamma^{k\ell,rs}_{ijpq}\neq 0$. Then, since $\ell\neq k$, $s\neq r$ and $(k,\ell)\neq (r,s)$, we deduce   from (\ref{gamfinal}) that we necessarily have $s=k$ and $\ell =r$. Thus $k<r$ and, since $r=\ell$, we obtain the inequality $k<\ell$. That is not possible since  $\ell\in\{1,\cdots,k-1\}$. So, we deduce that $\gamma^{k\ell,rs}_{ijpq}= 0$. Using (\ref{cov}), (\ref{gamma}) and (\ref{gamdiag}) we can conclude that $\Gamma =4h^{\prime\prime}(0)\,\mathbb{I}_{d}$. Then,  we deduce from Theorem 4.1 that $\mathcal{Q}=4h^{\prime\prime}(0)\mathcal{Q}^\prime$, where $\mathcal{Q}^\prime$ is a random variable with distribution equal to $\chi^2_d$.
\hfill $\Box$

\end{document}